%% file: tcx.tex
\documentclass[11pt,final]{article}%
\usepackage{stix2}
\usepackage{appendix}
\usepackage{amsfonts}
\usepackage{amsmath}
\usepackage{amssymb}
\usepackage{amsxtra,amscd}
\usepackage[amsmath,hyperref,thmmarks]{ntheorem}
\usepackage{makeidx}
\usepackage{graphicx}
\usepackage{cprotect}
\usepackage[colorlinks=true,bookmarksnumbered=true,pdfpagemode=None,final]%
{hyperref}%
\providecommand{\U}[1]{\protect\rule{.1in}{.1in}}
\theoremnumbering{arabic}
\theoremheaderfont{\scshape}
\RequirePackage{latexsym}
\theorembodyfont{\slshape}
\theoremseparator{.}

\newtheorem{theoreml}{Theorem}

\theorembodyfont{\upshape}

\theorembodyfont{\small}

\newtheorem*{daside}{Aside}

\theorembodyfont{\normalsize}
\theoremstyle{nonumberplain}
\theoremsymbol{\ensuremath{_\Box}}

\qedsymbol{\ensuremath{_\Box}}
\input{defa}

\begin{document}

\title{Tannakian categories: origins and summary}
\date{\today}
\author{J.S. Milne}
\maketitle

\section{Origins}

Andr\'{e} Weil's work on the arithmetic of curves and other varieties over
finite fields led him in 1949 to state his famous \textquotedblleft Weil
conjectures\textquotedblright. These had a profound influence on algebraic
geometry and number theory in the following decades. In an effort to explain
the conjectures, Grothendieck was led to define several different
\textquotedblleft Weil cohomology theories\textquotedblright\ and to posit an
ur-theory underlying all of them whose objects he called motives. In order to
provide a framework for studying these different theories, especially motives,
Grothendieck introduced the notion of a tannakian category.

Weil's first insight was that the numbers of points on smooth projective
algebraic varieties over finite fields behave as if they were the alternating sums of 
the traces of an
operator acting on a well-behaved homology theory.\footnote{Il me fallut du
temps avant de pouvoir m\^eme imaginer que les nombres de Betti fussent
susceptibles d'une interpr\'etation en g\'eom\'etrie alg\'ebrique abstraite.
Je crois que je fis un raisonnement heuristique bas\'e sur la formule de
Lefschetz. (It took me a while before I could even imagine that the Betti
numbers were susceptible to an interpretation in abstract algebraic geometry.
I think I made a heuristic argument based on the Lefschetz formula). Weil,
\OE uvre, Commentaire [1949b].}
In particular, the (co)homology groups should be vector spaces over a field of
characteristic zero, be functorial, and give the \textquotedblleft
correct\textquotedblright\ Betti numbers. However, already in the 1930s,
Deuring and Hasse had shown that the endomorphism algebra of an elliptic curve
over a field of characteristic $p$ may be a quaternion algebra over
$\mathbb{Q}$ that remains a division algebra even when tensored with
$\mathbb{Q}_{p}$ or $\mathbb{R}{}$, and hence cannot act on a $2$-dimensional
vector space over $\mathbb{Q}{}$ (or even $\mathbb{Q}{}_{p}$ or $\mathbb{R}$).
In particular, no such cohomology theory with $\mathbb{Q}{}$-coefficients exists.

Grothendieck defined \'{e}tale cohomology groups with $\mathbb{Q}{}_{\ell}%
$-coefficients for each prime $\ell$ distinct from the characteristic of the
ground field, and in characteristic $p\neq0$, he defined the crystalline
cohomology groups with coefficients in an extension of $\mathbb{Q}{}_{p}$.
Each cohomology theory is well-behaved. In particular it has a Lefschetz trace
formula, and Weil's first insight is explained by realizing the points of the
variety in a finite field as the fixed points of the Frobenius operator, and
hence, by the trace formula, their cardinality as the alternating sum of the 
traces of the Frobenius operator
acting on the cohomology groups. A striking feature of this is that, while the
traces of the Frobenius operator are, by definition, elements of different
fields $\mathbb{Q}{}_{l}$, they in fact lie in $\mathbb{Q}{}$ and are
independent of $l$ (for smooth projective varieties). This last fact suggested
to Grothendieck that there was some sort of $\mathbb{Q}{}$-theory underlying
the different $\mathbb{Q}_{l}$-theories. To explain what this is, we need the
notion of a tannakian category.

Briefly, a tannakian category over a field $k$ is a $k$-linear abelian
category with a tensor product structure having most of the properties of the
category of finite-dimensional representations of an affine group scheme over
$k$ except one: there need not exist an exact tensor functor to the category
of $k$-vector spaces, and when one does exist there need not be a canonical one.
Each of the cohomology theories takes values, not just in a category of vector
spaces, but in a tannakian category. For example, crystalline cohomology takes
values in a category of isocrystals. These are finite-dimensional vector
spaces over an extension of $\mathbb{Q}_{p}$, but only the elements of 
$\mathbb{Q}_{p}$ act as endomorphisms in the category. More specifically, if $\1$ is the unit object of the
category (the tensor product of the empty set), we have $\mathrm{End}(\1)=\mathbb{Q}%
_{p}$. Grothendieck's insight is that there should be a tannakian category
$\mathsf{Mot}$ over $\mathbb{Q}{}$ such that the functors to the local tannakian
categories defined by the different cohomology theories factor through it.
Algebraic correspondences between smooth projective algebraic varieties should
define maps between motives, whose traces lie in $\mathrm{End}(\1)=\mathbb{Q}$ and map to the
traces on the various cohomology groups, which explains why the latter lie in
$\mathbb{Q}$.

Weil's second insight was that an analogue of the Riemann hypothesis should
hold for the eigenvalues of Frobenius operators. This suggested that some of
the well-known positivities in characteristic zero should persist to
characteristic $p$. To see why, we briefly recall Weil's proof of the Riemann
hypothesis for abelian varieties over finite fields.

Consider an abelian variety $A$ over an algebraically closed field of
characteristic $p$. For a prime $\ell\neq p$, we have a finite-dimensional
$\mathbb{Q}{}_{\ell}$-vector space $V_{\ell}A$, and, for each polarization of
$A$, we have a pairing $\varphi\colon V_{\ell}A\times V_{\ell}A\rightarrow
\mathbb{Q}{}_{\ell}$. As $\mathbb{Q}{}_{\ell}$ is not a subfield of
$\mathbb{R}{}$, it makes no sense to say that $\varphi$ is positive-definite.
However, Weil showed that $\varphi$ induces an involution on the
finite-dimensional $\mathbb{Q}{}$-algebra $\mathrm{End}(A)\otimes\mathbb{Q}{}$ and
that this involution \textit{is} positive.\footnote{Over $\mathbb{C}$, this
was known to the Italian geometers as the positivity of the Rosati
involution.} The Riemann hypothesis for the abelian variety follows directly
from this. Grothendieck extended Weil's ideas to tannakian categories by
introducing the notion of a ``Weil form'' on an object of a tannakian category
and of a ``polarization'' on a tannakian category.

A tannakian category over $k$ is said to be neutral if it admits an exact
tensor functor to the category of $k$-vector spaces. Neutral tannakian
categories are the analogues for affine group schemes of the categories
studied by Tannaka and Krein. A classical theorem of Tannaka describes how to
recover a compact topological group from its category of finite-dimensional
unitary representations, and Krein characterized the categories arising in
this way.

Not all tannakian categories are neutral, and the obstruction to a tannakian
category over $k$ having a $k$-valued fibre functor lies in a nonabelian
cohomology group of degree 2, more general than was available in the early
1960s. Grothendieck's student Giraud developed the necessary nonabelian
cohomology theory in his thesis.\footnote{Giraud, Jean, 
Cohomologie non ab\'elienne. 
Die Grundlehren der mathematischen Wissenschaften, Band 179. Springer-Verlag, 
Berlin-New York, 1971.}

As we have explained, the idea of tannakian categories, and of their
importance for motives, was Grothendieck's. He explained it to Saavedra
Rivano, who developed the theory of tannakian categories in his 
thesis.\footnote{Saavedra Rivano, Neantro. Cat\'egories 
Tannakiennes. Lecture Notes in Mathematics, Vol. 265. Springer-Verlag, 
Berlin-New York, 1972.} It was Saavedra who introduced the terminology
\textquotedblleft tannakian\textquotedblright. Although Grothendieck used the
term \textquotedblleft tannakian category\textquotedblright\ in unpublished
writings, he considered the categories to be part of a vast theory
engobalizing Galois theory and the theory of fundamental groups, and later
wrote that \textquotedblleft Galois--Poincar\'{e} category\textquotedblright%
\ would have been a more appropriate name.\footnote{Deligne writes: I expect
that at first Grothendieck did not know of Tannaka's work -- and never cared
about it. His aim was to unify the cohomology theories he had created. That
each $H$ is with values in a category with $\otimes$, and that K\"unneth
holds, was a brilliant insight which, like a number of his brilliant ideas, is
now part of our subconscious, making it hard to see how deep it was.}

\section{Summary}

We now present a summary of the main results of the theory. Throughout, $k$ is
a field.

A \emph{tensor category}%
\index{tensor category}
(symmetric monoidal category) is a category $\mathsf{C}$ together with a
functor $\otimes\colon\mathsf{C}{}\times\mathsf{C}\rightarrow\mathsf{C}$ and
sufficient constraints to ensure that the tensor product of any (unordered)
finite set of objects in $\mathsf{C}$ is well-defined up to a canonical
isomorphism. In particular, there exists a unit object $\1$ (tensor product of
the empty set of objects). A tensor category is \emph{rigid}%
\index{rigid}
if every object admits a dual (in a strong sense). A \emph{tensor functor}%
\index{tensor functor}
of tensor categories is one preserving the tensor products and constraints.

A \emph{tensorial category over }$k$%
\index{tensorial category}
is a rigid abelian tensor category equipped with a $k$-linear structure such
that $\otimes$ is $k$-bilinear and the structure map $k\to\mathrm{End}(\1)$ is an
isomorphism. A tensorial category over $k$ is a \emph{tannakian category over
}$k$%
\index{tannakian category}
if, for some nonzero $k$-algebra $R$, there exists an $R$\emph{-valued}
\emph{fibre functor}%
\index{fibre functor}%
, i.e., an exact $k$-linear tensor functor $\omega\colon\mathsf{C}%
\rightarrow\mathrm{Mod}(R)$. We write $\mathrm{Aut}^{\otimes}(\omega)$ for the group of
automorphisms of $\omega$ (as a tensor functor).

In the remainder of the introduction, all tensor categories are assumed to be
essentially small (i.e., equivalent to a small category).

\subsection{A criterion to be a tannakian category}

For an object $X$ of a tensorial category $\mathsf{C}$ over $k$, there is a
canonical trace map
\[
\mathrm{Tr}_{X}\colon\mathrm{End}(X)\rightarrow\mathrm{End}(\1)\simeq k,
\]
and we let $\dim X$ denote the trace of $\mathrm{id}_{X}$. In tensorial categories,
traces are additive on short exact sequences (I, 6.6).
\footnote{All references in this section are to: Milne, J.\,S., Tannakian
Categories (version February 16, 2025), available 
\href{https://jmilne.org/math/Books/tcdraft.pdf}{here}, from which
this article has been adapted.}

\begin{theoreml}
[I, 10.1]\label{T5}A tensorial category over $k$ of characteristic zero
is tannakian (i.e., a fibre functor exists) if and only if, for all objects
$X$, $\dim X$ is an integer $\geq0$.
\end{theoreml}

\subsection{Neutral tannakian categories.}

A tannakian category $(\mathsf{C},\otimes)$ over $k$ is \emph{neutral}%
\index{tannakian category!neutral}
if there exists a $k$-valued fibre functor. For example, the category
$\mathsf{Repf}(G)$ of finite-dimensional representations of an affine group scheme $G$
over $k$ is a tannakian category over $k$ with the forgetful functor as a
$k$-valued fibre functor.

\begin{theoreml}
[II, 3.1]\label{T1} Let $\mathsf{C}$ be a tannakian category over $k$
and $\omega$ a $k$-valued fibre functor.

\begin{enumerate}
\item The functor of $k$-algebras $R\rightsquigarrow\mathrm{Aut}^{\otimes}%
(\omega\otimes R)$ is represented by an affine group scheme $G=\mathcal{A}%
ut^{\otimes}_k(\omega)$ over $k$.

\item The functor $\mathsf{C}\rightarrow\mathsf{Repf}(G)$ defined by $\omega$ is an
equivalence of tensor categories.
\end{enumerate}
\end{theoreml}

\noindent For example, if $\mathsf{C}=\mathsf{Repf}(G)$ and $\omega$ is the forgetful
functor, then $\mathcal{A}ut^{\otimes}_{k}(\omega)=G$.

The theorem gives a dictionary between neutralized tannakian categories over
$k$ and affine group schemes over $k$. To complete the theory in the neutral
case, it remains to describe the $R$-valued fibre functors on $\mathsf{C}$ for
$R$ a $k$-algebra.

\begin{theoreml}
[II, 8.1]\label{T2} Let $\mathsf{C}$ and $\omega$ be as in Theorem
\ref{T1}, and let $G=\mathcal{A}ut^{\otimes}_{k}(\omega)$. For any $R$-valued
fibre functor $\nu$ on $\mathsf{C}$, $\mathcal{I}som^{\otimes}(\omega\otimes
R,\nu)$ is a torsor under $G_{R}$ for the fpqc topology. The functor
$\nu\rightsquigarrow\mathcal{I}som^{\otimes}(\omega\otimes R,\nu)$ is an
equivalence from the category of $R$-valued fibre functors on $\mathsf{C}$ to
the category of $G_{R}$-torsors,%
\[
\Fib(\mathsf{C})_{R}\sim\Tors(G)_{R}.
\]

\end{theoreml}

\begin{daside}
\label{d3.4}The situation described in the theorem is analogous to the
following. Let $X$ be a connected topological space, and let $\mathsf{C}$ be
the category of locally constant sheaves of $\mathbb{Q}{}$-vector spaces on
$X$. For each $x\in X$, there is a functor $\omega_{x}\colon\mathsf{C}%
\rightarrow\mathsf{Vecf}_{\mathbb{Q}{}}$ taking a sheaf to its fibre at $x$, and
$\omega_{x}$ defines an equivalence of categories $\mathsf{C}\rightarrow
\mathsf{Rep}_{\mathbb{Q}}(\pi_{1}(X,x))$. Let $\Pi_{x,y}$ be the set of homotopy
classes of paths from $x$ to $y$; then $\Pi_{x,y}\simeq\mathrm{Isom}(\omega_{x}%
,\omega_{y})$, and $\Pi_{x,y}$ is a $\pi_{1}(X,x)$-torsor.
\end{daside}

\subsection{General tannakian categories.}

Many of the tannakian categories arising in algebraic geometry are not
neutral. They correspond to affine \textit{groupoid} schemes rather than
affine \textit{group} schemes.

Let $S$ be an affine scheme over $k$. A $k$-\emph{groupoid scheme acting on}%
\index{groupoid scheme}
$S$ is a $k$-scheme $G$ together with two $k$-morphisms $t,s\colon
G\rightrightarrows S$ and a partial law of composition
\[
\circ\colon G\underset{s,S,t}{\times}G\rightarrow G\quad\quad\text{(morphism
of S}\times_{k}\text{S-schemes){}}%
\]
such that, for all $k$-schemes $T$, $(S(T),G(T),(t,s),\circ)$ is a groupoid
(i.e., a small category in which the morphisms are isomorphisms). The groupoid
$G$ is \emph{transitive}%
\index{groupoid!transitive}
if the morphism
\[
(t,s)\colon G\rightarrow S\times_{k}S
\]
is faithfully flat. The representations of $G$ on locally free sheaves of
finite rank on $S$ form a tannakian category $\mathsf{Repf}(S\mathrm{:}G)$ over $k$.

Let $S=\mathrm{Spec}(R)$ be an affine scheme over $k$. By a fibre functor over $S$, we
mean an $R$-valued fibre functor. For example, $\mathsf{Repf}(S\mathrm{:}G)$ has a
canonical (forgetful) fibre functor over $S$. When $\omega$ is a fibre functor
over $S$ on a tannakian category over $k$, we let $\mathcal{A}ut_{k}^{\otimes
}(\omega)$ denote the functor of $S\times_{k}S$-schemes sending $(b,a)\colon
T\rightarrow S\times_{k}S$ to $\mathrm{Isom}_{T}^{\otimes}(a^{\ast}\omega,b^{\ast
}\omega)$.

\begin{theoreml}
[III, 1.1]\label{T3} Let $\mathsf{C}$ be a tannakian category over $k$
and $\omega$ a fibre functor over $S$.

\begin{enumerate}
\item The functor $\mathcal{A}ut_{k}^{\otimes}(\omega)$ is represented by an
affine $k$-groupoid scheme $G$ acting transitively on $S$.

\item The functor $\mathsf{C}\rightarrow\mathrm{Repf}(S\mathrm{:}G)$ defined by
$\omega$ is an equivalence of tensor categories.
\end{enumerate}
\end{theoreml}

\noindent For example, if $\mathsf{C}=\mathrm{Repf}(S\mathrm{:}G)$ and $\omega$ is the
forgetful functor, then $\mathcal{A}ut^{\otimes}_{k}(\omega)\simeq G$.

\subsection{The gerbe of fibre functors}

Let $\mathrm{Aff}_{k}$ denote the category of affine $k$-schemes. For each affine
$k$-scheme $S$, we let $\Fib(\mathsf{C})_{S}$ denote the category of fibre
functors on $\mathsf{C}$ over $S$. As $S$ varies, the categories
$\Fib(\mathsf{C})_{S}$ form a stack over $\mathrm{Aff}_{k}$ for the fpqc topology, and
(c) of Theorem \ref{T3} implies that $\Fib(\mathsf{C})$ is a gerbe (any two
fibre functors are locally isomorphic).

The tannakian categories over $k$ form a $2$-category with the $1$-morphisms
being the exact $k$-linear tensor functors and the $2$-morphisms the morphisms
of tensor functors. Similarly, the affine gerbes\footnote{A gerbe
is affine if the automorphisms of any object form an affine group scheme.} over $k$ form a $2$-category
with the $1$-morphisms being the cartesian functors of fibred categories and
the $2$-morphisms being the equivalences between $1$-morphisms.

\begin{theoreml}
[IV, 3.3]\label{T4} The $2$-functor sending a tannakian category to its
gerbe of fibre functors is an equivalence of $2$-categories.\footnote{Not a
$2$-equivalence} Explicitly, for any tannakian category $\mathsf{C}$ over $k$,
the canonical functor%
\[
\mathsf{C}\rightarrow\mathrm{Repf}(\Fib(\mathsf{C}))
\]
is an equivalence of tensor categories, and for any affine gerbe $\mathsf{G}$
over $k$, the canonical functor%
\[
\mathsf{G}\rightarrow\Fib(\mathrm{Repf}(\mathsf{G}))
\]
is an equivalence of stacks.
\end{theoreml}

The theorem gives a dictionary between tannakian categories over $k$ and
affine gerbes over $k$.

\subsection{The fundamental group of a tannakian category}

Let $\mathsf{T}$ be a tannakian category over $k$. The notion of a Hopf
algebra makes sense in the ind-category $\mathrm{Ind}\mathsf{T}$. In order to make
available a geometric language, Deligne defined the category of affine group
schemes in $\mathrm{Ind}\mathsf{T}{}$ to be the opposite of that of commutative Hopf
algebras. If $G$ is the group scheme corresponding to the Hopf algebra $A$,
then, for any $R$-valued fibre functor $\omega$, $\omega(G)\overset{\df}{=}%
\mathrm{Spec}(\omega(A))$ is an affine group scheme over $R$. The \emph{fundamental
group}%
\index{fundamental group}
$\pi(\mathsf{T})$ of $\mathsf{T}$ is the affine group scheme in
$\mathrm{Ind}\mathsf{T}$ such that%
\[
\omega(\pi(\mathsf{T}))=\mathcal{A}ut^{\otimes}(\omega)
\]
for all fibre functors $\omega$. The group $\pi(\mathsf{T})$ acts on the
objects $X$ of $\mathsf{T}$, and $\omega$ transforms this action into the
natural action of $\mathcal{A}ut^{\otimes}(\omega)$ on $\omega(X)$.

Let $X$ be a topological space, connected, locally connected, and locally
simply connected. There is the following analogy:\smallskip

\begin{center}%
\begin{tabular}
[c]{l|l}%
$\qquad\mathsf{T}$ & $\qquad X$\\\hline
object $Y$ of $\mathsf{T}$ & covering of $X$(=locally constant sheaf)\\\hline
fibre functor $\omega_{0}$ & point $x_{0}\in X$\\\hline
$\mathcal{A}ut^{\otimes}(\omega_{0})$ & $\pi_{1}(X,x_{0})$\\\hline
$\pi(\mathsf{T})$ & local system of the $\pi_{1}(X,x)$\\\hline
action of $\pi(\mathsf{T})$ on $Y$ in $\mathsf{T}$ & action of the local
system of the $\pi_{1}(X,x)$\\
& on a locally constant sheaf.
\end{tabular}

\end{center}

\smallskip\noindent For $\mathsf{T}$ the category of motives over $k$,
$\pi(\mathsf{T})$ is called the \emph{motivic Galois group}%
\index{motivic Galois group}
of $k$.\footnote{From Deligne: The first three lines [in the table] were
surely clear and important for Grothendieck. I don't remember him considering
$\mathrm{Ind}\mathsf{T}$, $\pi(\mathsf{T})$, or Hopf algebras in $\mathsf{T}$. For me,
it was a way to make sense of my surprise, seeing that for each of the
standard fibre functors $\omega$ with values in $\mathcal{C}$,
\[
\mathcal{A}ut^{\otimes}(\omega\colon\text{motives}\to\mathcal{C}%
\to\text{vector spaces})
\]
had the same `texture' as objects of $\mathcal{C}$.}

\subsection{Polarized tannakian categories.}

For tannakian categories over $\mathbb{R}{}$ (or a subfield of $\mathbb{R}{}%
$), there are positivity structures called \emph{polarizations}.%
\index{polarization}
For simplicity, let $(\mathsf{C},\otimes)$ be an algebraic tannakian category
over $\mathbb{R}$. A nondegenerate bilinear form%
\[
\phi\colon V\otimes V\rightarrow\mathbb{R}{}%
\]
on an object $V$ of $\mathbb{C}{}$ is called a \emph{Weil form}%
\index{Weil form}
if its parity $\epsilon_{\phi}$ (the unique automorphism of $V$ satisfying
$\phi(y,x)=\phi(x,\epsilon_{\phi}y)$) is in the centre of $\mathrm{End}(V)$ and if for
all nonzero endomorphisms $u$ of $V$, $\mathrm{Tr}(u\circ u^{\phi})>0$, where
$u^{\phi}$ is the adjoint of $u$. Two Weil forms $\phi\colon V\otimes
V\rightarrow\mathbb{R}{}$ and $\psi\colon W\otimes W\rightarrow\mathbb{R}{}$
are \emph{compatible}%
\index{compatible Weil forms}
if the form $\phi\oplus\psi$ on $V\oplus W$ is again a Weil form.

Now fix an $\epsilon\in Z(\mathbb{R}{})$, where $Z$ is the centre of the band
of the gerb of $\Fib(\mathsf{C})$ -- it is a commutative algebraic
$\mathbb{R}{}$-group -- and suppose that for each object $V$ of $\mathsf{C} $
we are given a nonempty compatibility class $\pi(V)$ of ($\pi$\emph{-positive}%
) Weil forms on $V$ with parity $\epsilon_{V}$. We say that $\pi$ is an
$\epsilon$\emph{-polarization} of $\mathsf{C}$ if direct sums and tensor
products of $\pi$-positive forms are $\pi$-positive. When $\epsilon=1$, so
that $\phi(x,y)=\phi(y,x)$, the polarization is said to be \emph{symmetric}.%
\index{polarization!symmetric}%

Let $G$ be an affine group scheme over $\mathbb{R}$, and let $C$ be an element
of $G(\mathbb{R})$ such that $\mathrm{inn}(C)$ is a Cartan involution, i.e., the
involution corresponding to a compact form\footnote{A real form $G^{\prime}$
of $G$ is \emph{compact}%
\index{compact real form}
if $G(\mathbb{R})$ is compact and contains a point of each connected component
of $G_{\mathbb{C}}$.} of $G$. Because $\mathrm{inn}(C)$ is an involution, $C^{2}$ is
central. For each $V$ in $\mathsf{Repf}(G)$, let $\pi_{C}(V)$ be the set of
$G$-invariant bilinear forms $\phi\colon V\otimes V\rightarrow\mathbb{R}{}$
such that the bilinear form $\phi_{C}$,%
\[
\phi_{C}(x,y)\overset{\df}{=}\phi(x,Cy)\text{,}%
\]
is symmetric and positive-definite. Then $\pi_{C}$ is a $C^{2}$-polarization
on $\mathsf{Repf}(G)$. For a neutralized tannakian category, the $\pi_{C}$ exhaust the polarizations.

\begin{theoreml}[V, 8.2]
\label{T6} Let $G$ be an affine algebraic $\mathbb{R}{}$-group. Then
$\mathsf{Repf}(G)$ admits a polarization if and only if $G$ is an inner form of a real
compact group, in which case every polarization is of the form $\pi_{C}$ for
some $C $ as above, and $C$ is uniquely determined by the polarization up to conjugacy.
\end{theoreml}

It follows from the theorem that if $\mathsf{C}$ is an algebraic tannakian
category endowed with a symmetric polarization, then $\mathsf{C}$ is neutral
and there is an $\mathbb{R}{}$-valued fibre functor $\omega\colon\mathsf{C}%
{}\rightarrow\mathsf{Vecf}(\mathbb{R}{})$ such that
$\mathcal{A}ut^{\otimes}(\omega)$ is a compact $\mathbb{R}$-group; moreover,
$\omega$ is uniqely determined up to a unique isomorphism by the condition
that the positive forms on an object $V$ of $\mathsf{C}$ are exactly the forms
$\phi$ such that $\omega(\phi)$ is symmetric and positive-definite.

\subsection{Motives}

Fix an admissible equivalence relation for algebraic cycles on smooth
projective algebraic varieties over $k$, and let $\mathsf{M}(k)$ denote the
corresponding category of motives. It is a tensor category equipped with a
$\mathbb{Q}$-linear structure (in particular, it is additive) such that
$\otimes$ is $\mathbb{Q}$-bilinear.

\begin{theoreml}[VI, 2.5]
\label{T7}The category of motives $\mathsf{M}(k)$ is a $\mathbb{Q}$-linear
rigid tensor category.
\end{theoreml}

Let $X$ be a smooth projective variety over $k$. We say that $X$ satisfies the
\emph{sign conjecture}%
\index{sign conjecture}
if there exists an algebraic cycle $e$ on $X\times X$ such that $e^2=e$ and $eH^{\ast
}(X)=\bigoplus_{i\geq0} H^{2i}(X)$ for the standard Weil cohomology theories.
Smooth projective varieties over a finite field satisfy the sign conjecture,
as do abelian varieties over any field. Let $\mathsf{NMot}(k)$ denote the
category of motives for numerical equivalence over $k$ generated by the smooth
projective varieties over $k$ satisfying the sign conjecture.

\begin{theoreml}[VI, 6.12]
\label{T8} The category of numerical motives $\mathsf{NMot}(k)$ is a
semisimple tannakian category over $\mathbb{Q}$.
\end{theoreml}

To prove that $\mathsf{NMot}(k)$ is polarized and that the standard Weil cohomologies
factor through it requires Grothendieck's standard conjectures. Given the lack
of progress on these conjectures, Deligne has suggested looking for
alternatives, of which there are several. 

\subsection{Acknowledgements}
I thank Pierre Deligne for his generous help.

\end{document}

%% file: defa.tex
\def\1{\mathbb{1}}
\newcommand{\Tors}{\textsc{Tors}}
\newcommand{\Fib}{\textsc{Fib}}

\newcommand{\df}{\smash{\lower.12em\hbox{\textup{\tiny def}}}}

\usepackage{changepage}

\let\quoteOld\quote
\let\endquoteOld\endquote
\renewenvironment{quote}{\small\quoteOld}{\endquoteOld}

\usepackage[overload]{textcase}

\usepackage{xcolor}
\definecolor{bblue}{rgb}{0.0, 0.0, 0.6}

\usepackage{array}

\usepackage{microtype}

\usepackage{tikz}
\usepackage{tikz-cd}
\usetikzlibrary{matrix,arrows,positioning,decorations.pathmorphing}
\usetikzlibrary{decorations.markings,shapes.geometric}
\tikzset{commutative diagrams/column sep/Huge/.initial=12ex}
\tikzcdset{arrow style=math font}

\usepackage[shortlabels]{enumitem}
\setitemize[1]{label=$\diamond$\hspace{0.07in}}
\setlist{topsep=0.1em,itemsep=0.1em,parsep=0.1em}

\setdescription{font=\normalfont}

\usepackage{url}
\usepackage{verbatim}

\usepackage[notcite,notref]{showkeys}

\usepackage{mathtools}

\newcommand{\eb}[1]{{\itshape\bfseries#1}}
\renewcommand{\emph}{\eb}

\usepackage{titletoc,titlesec}
\contentsmargin{2.0em}
\dottedcontents{section}[2.3em]{}{1.8em}{1pc}
\usepackage[nottoc]{tocbibind}

\titleformat*{\section}{\LARGE\bfseries}
\titleformat*{\subsection}{\Large\itshape}
\titleformat*{\subsubsection}{\scshape}
\titleformat*{\paragraph}{\itshape}

\usepackage{natbib}
\let\cite\citealt

\hypersetup{linkcolor=bblue,anchorcolor=bblue,citecolor=bblue,filecolor=bblue,urlcolor=bblue}
\hypersetup{pdftitle={Shimura Varieties},pdfauthor={J.S. Milne}}
\hypersetup{pdfsubject={none},pdfkeywords={none}}

\usepackage[papersize={6.5in,10.0in},margin=0.5in,pdftex]{geometry}